# Spectral Differentiation Operators And Hydrodynamic Models For Stability Of Swirling Fluid Systems

DIANA ALINA BISTRIAN
Department of Electrical Engineering and Industrial Informatics
Engineering Faculty of Hunedoara,"Politehnica" University of Timisoara
Str. Revolutiei Nr.5, Hunedoara, 331128
ROMANIA
diana.bistrian@fih.upt.ro

FLORICA IOANA DRAGOMIRESCU
Department of Mathematics
"Politehnica" University of Timisoara
Victoriei Square Nr.2, Timisoara, 300006
ROMANIA
ioana.dragomirescu@mat.upt.ro

GEORGE SAVII
Department of Mechatronics
Mechanical Engineering Faculty, "Politehnica" University of Timisoara
Mihai Viteazu Nr.1, Timisoara, 300222
ROMANIA
george.savii@mec.upt.ro

*Abstract:* In this paper we develop hydrodynamic models using spectral differential operators to investigate the spatial stability of swirling fluid systems. Including viscosity as a valid parameter of the fluid, the hydrodynamic model is derived using a nodal Lagrangean basis and the polynomial eigenvalue problem describing the viscous spatial stability is reduced to a generalized eigenvalue problem using the companion vector method. For inviscid study the hydrodynamic model is obtained by means of a class of shifted orthogonal expansion functions and the spectral differentiation matrix is derived to approximate the discrete derivatives. The models were applied to a Q-vortex structure, both schemes providing good results.

## 1 Introduction

The role of the hydrodynamic stability theory in fluid mechanics reach a special attention, especially when reaserchers deal with problem of minimum consumption of energy. This theory deserves special mention in many engineering fields, such as the aerodynamics of profiles in supersonic regime, the construction of automation elements by fluid jets and the technique of emulsions.

The main interest in recent decades is to use the theory of hydrodynamic stability in predicting transitions between laminar and turbulent configurations for a given flow field. R.E. Langer [1] proposed a theoretical model for transition based on supercritical branching of the solutions of the Navier-Stokes equations. This model was substantiated mathematically by E. Hopf [2] for systems of nonlinear equations close to Navier-Stokes equations. C.C. Lin, a famous specialist in hydrodynamic stability theory, published his first paper on stability of fluid systems in which the mathematical formulation of the problems was essentially diferent from the conservative treatment [3]. The Nobel laureate Chandrasekhar [4] presents in his study considerations of typical problems in hydrodynamic and hydromagnetic stability as a branch of experimental physics. Among the subjects treated are thermal instability of a layer of fluid heated from below, the Benard problem, stability of Couette flow, and the Kelvin-Helmholtz instability.

Many publications in the field of hydrodynamics are focused on vortex motion as one of the basic states of a flowing continuum and effects that vortex can produce. Mayer [5] and Korrami [6] have mapped out the stability of Q-vortices, identifying both inviscid and viscous modes of instability. The mathematical description of the dynamics of swirling flows is hindered by the requirement to consider three-dimensional and nonlinear effects, singularity and various instabilities as in [7, 8, 9].

The objective of this paper is to present new instruments that can provide relevant conclusions on the stability of swirling flows, assessing both analytically methodology and numerical methods. The study involves new mathematical models and simulation algorithms that translate equations into computer code instructions immediately following problem formulations. Classical vortex problems were chosen to validate the code with the existing results in the literature. The paper is outlined as follows: Section 1 gives a brief motivation for the study of hydrodynamic stability using computer aided techniques. The dispersion equation governing the linear stability analysis for swirling flows against normal mode perturbations is derived in Section 2. In Section 3 a nodal collocation method is proposed for viscous stability investigations and in Section 4 a modal collocation method is developed, based on shifted orthogonal expansions, assessing different boundary conditions. In Section 5 the hydrodynamic models are applied upon the velocity profile of a Q-vortex and Section 6 concludes the paper.

## 2 Problem Formulation

Hydrodynamic stability theory is concerned with the response of a laminar flow to a disturbance of small or moderate amplitude. If the flow returns to its original laminar state one defines the flow as stable, whereas if the disturbance grows and causes the laminar flow to change into a different state, one defines the flow as unstable. Instabilities often result in turbulent fluid motion, but they may also take the flow into a different laminar, usually more complicated state. The equations governing the general evolution of fluid flow are known as the Navier-Stokes equations [4]. They describe the conservation of mass and momentum. The evolution equations for the disturbance can be derived by considering a basic state $\{\overline{U} = (u_z, u_r, u_\theta), p\}$ and a perturbed state $\{\overline{V} = (v_z, v_r, v_\theta), \pi\}$, with the disturbance being of order $0 \prec \delta \prec\prec 1$

$$\overline{U} = [U(r), 0, W(r), P(r)] + \delta \overline{V} \quad (1)$$

Consistent with the parallel mean flow assumption is that the functional form for the mean part of the velocity components only involves the cross-stream coordinate and also zero mean radial velocity. The linearized equations are obtained after substituting the expressions for the components of the velocity and pressure field into the Navier Stokes equations and only considering contributions of first order in delta. For high Reynolds numbers a restrictive hypothesis to neglect viscosity can be imposed in some problems. The linearized equations in operator form are

$$L \cdot S = 0, \quad S = (v_r \quad v_\theta \quad v_z \quad \pi)^T \quad (2)$$

and the elements of matrix $L$ being

$L_{11} = \partial_t + \frac{W}{r}\partial_\theta + U\partial_z - \frac{1}{\text{Re}}\Delta - \frac{1}{r^2\text{Re}}$

$L_{12} = -\frac{2W}{r} + \frac{2}{r^2\text{Re}}\partial_\theta$, $L_{13} = 0$, $L_{14} = \partial_r$,

$L_{21} = W' + \frac{W}{r} - \frac{2}{r^2\text{Re}}\partial_\theta$,

$L_{22} = \partial_t + \frac{W}{r}\partial_\theta + U\partial_z - \frac{1}{\text{Re}}\Delta + \frac{1}{r^2\text{Re}}$, $L_{23} = 0$,

$L_{24} = \frac{1}{r}\partial_\theta$, $L_{31} = U'$,

$L_{32} = 0$, $L_{33} = \partial_t + \frac{W}{r}\partial_\theta + U\partial_z - \frac{1}{\text{Re}}\Delta$,

$L_{34} = \partial_z$, $L_{41} = \partial_r + \frac{1}{r}$,

$L_{42} = \frac{1}{r}\partial_\theta$, $L_{43} = \partial_z$, $L_{44} = 0$,

where $\partial_{\{t,z,r,\theta\}}$ denote the partial derivative operators and primes denote derivative with respect to radial coordinate. In linear stability analysis the disturbance components of velocity are shaped into normal mode form, given here

$$\{v_z, v_r, v_\theta, \pi\} = \{F(r), iG(r), H(r), P(r)\}E(t,z,\theta) \quad (3)$$

where $E(t,z,\theta) \equiv e^{i(kz+m\theta-\omega t)}$, $F,G,H,P$ represent the complex amplitudes of the perturbations, $k$ is the complex axial wave number, $m$ is the tangential integer wave number and $\omega$ represents the complex frequency. The hydrodynamic equation of dispersion is obtained, where we have explicitly decomposed into operators that multiply $\omega$ and the different powers of $k$

$$(\omega M_\omega + M + kM_k + k^2 M_{k^2}) \cdot (F \quad G \quad H \quad P)^T = 0. \quad (4)$$

The non zero elements of matrices are given explicitly by

$M_{\omega 1,1} = 1$, $M_{\omega 2,2} = -i$, $M_{\omega 3,3} = -i$,

$M_{k^2 1,1} = i/\text{Re}$, $M_{k^2 2,2} = i/\text{Re}$, $M_{k^2 3,3} = i/\text{Re}$,

$M_{k 1,1} = -U$, $M_{k 2,2} = Ui$, $M_{k 3,3} = Ui$,

$M_{k 3,4} = i$, $M_{k 4,3} = i$

and the elements of matrix $M$ are

$$M_{11} = -\frac{mW}{r} - \frac{i}{\text{Re}}d_{rr} - \frac{i}{r\text{Re}}d_r + \frac{i(m^2+1)}{r^2\text{Re}},$$

$$M_{12} = -\frac{2W}{r} + \frac{2im}{r^2\text{Re}}, \; M_{13} = 0, \; M_{14} = d_r,$$

$$M_{21} = iW' + \frac{iW}{r} + \frac{2m}{r^2\text{Re}},$$

$$M_{22} = \frac{imW}{r} - \frac{1}{\text{Re}}d_{rr} - \frac{1}{r\text{Re}}d_r + \frac{m^2}{r^2\text{Re}}, \; M_{23} = 0,$$

$$M_{24} = \frac{im}{r}, \; M_{31} = iW', \; M_{32} = 0,$$

$$M_{33} = \frac{imW}{r} - \frac{1}{\text{Re}}d_{rr} - \frac{1}{r\text{Re}}d_r + \frac{m^2}{r^2\text{Re}},$$

$$M_{34} = 0, \; M_{41} = id_r + \frac{i}{r}, \; M_{42} = \frac{im}{r}, \; M_{43} = M_{44} = 0,$$

where prime denotes differentiation with respect to the radius and $d_r$ and $d_{rr}$ mean the differentiation operators of first and second order.

## 3 Nodal Collocation Approach For Spatial Stability Including Viscosity

When the complex frequency $\omega = \omega_r + i \cdot \omega_i$, $\omega_r = \text{Re}(\omega)$, $\omega_i = \text{Im}(\omega)$ is determined as a function of the real wave number $k$ a temporal stability analysis is performed. Conversely, solving the dispersion relation (4) for the complex wave number $k = k_r + i \cdot k_i$, $k_r = \text{Re}(k)$, $k_i = \text{Im}(k)$, when $\omega$ is given real leads to the spatial branches $k(\omega, \Upsilon)$ where by $\Upsilon$ we denoted the set of all other physical parameters involved. In both cases, the sign of the imaginary part indicates the decay or either the growth of the disturbance. The growth of the wave solution in spatial case depends on the imaginary part of the axial wavenumber, as described in the next formula

$$e^{-k_iz}\left\{\begin{array}{l}F_r\cos(k_rz+\Theta)-F_i\sin(k_rz+\Theta)+\\i[F_r\sin(k_rz+\Theta)+F_i\cos(k_rz+\Theta)]\end{array}\right\}, \; \Theta \equiv m\theta-\omega t. \;(5)$$

Here the flow is considered unstable when the disturbance grows, i.e. the imaginary part of $k$ is negative. A given $\omega$ leads to a polynomial eigenvalue problem of form

$$\left(M_0 + kM_k + k^2M_{k^2}\right)\left(F \; G \; H \; P\right)^T = 0, M_0 \equiv \omega M_\omega + M. \;(6)$$

In general, the direct solution of polynomial eigenvalue problems can be heavy. For this case, we can transform the polynomial eigenvalue problem into a generalized eigenvalue problem, using the companion vector method, assessed also in [10]. We augment the system with the variable

$$\widehat{\Psi} \equiv \left(kF \; kG \; kH\right)^T, \; \widehat{S} \equiv \left(F \; G \; H \; P\right)^T \quad (7)$$

The eigenvalue problem describing the spatial hydrodynamic stability for a viscous fluid system reads now

$$\begin{pmatrix}M_0 & 0\\0 & I\end{pmatrix}\begin{pmatrix}\widehat{S}\\\widehat{\Psi}\end{pmatrix} + k\begin{pmatrix}M_k & M_{k^2}\\-I & 0\end{pmatrix}\begin{pmatrix}\widehat{S}\\\widehat{\Psi}\end{pmatrix} = 0 \quad (8)$$

where the first row is the polynomial eigenvalue problem (6) and the second row enforces the definition of $\widehat{\Psi}$.

The collocation method is associated with a grid of clustered nodes $x_j$ and weights $w_j$ $(j=0,...,N)$. The collocation nodes must cluster near the boundaries to diminish the negative effects of the Runge phenomenon [11]. Another aspect is that the convergence of the interpolation function on the clustered grid towards unknown solution is extremely fast. We recall that the nodes $x_0$ and $x_N$ coincide with the endpoints of the interval $[a,b]$, and that the quadrature formula is exact for all polynomials of degree $\leq 2N-1$, i. e.,

$$\sum_{j=0}^{N}v(x_j)w_j = \int_a^b v(x)w(x)dx, \quad (9)$$

for all $v$ from the space of test functions.

Let $\{\Phi_\ell\}_{\ell=0..N}$ a finite basis of polynomials relative to the given set of nodes, not necessary being orthogonal. If we choose a basis of non-orthogonal polynomials we refer to it as a *nodal* basis, Lagrange polynomials for example. In nodal approach, each function of the nodal basis is responsible for reproducing the value of the polynomial at one particular node in the interval. When doing simulations and solving PDEs, a major problem is one of representing a deriving functions on a computer, which deals only with finite integers. In order to compute the radial and pressure derivatives that appear in our mathematical model, the derivatives are approximated by differentiating a global interpolative function built trough the collocation points. We choose $\{\Phi_i\}_{i=0..N}$ given by Lagrange's formula

$$\Phi_i(r) = \frac{\omega_N(r)}{\omega'_N(r_i)(r-r_i)}, \text{ where } \omega_{N(r)} = \prod_{m=1}^{N}(r-r_m).$$

We constructed the interpolative spectral differentiation matrix $\Delta_{(N+1)\times(N+1)}$, having the entries

$$\Delta_{00} = \frac{2N^2+1}{6}, \Delta_{NN} = -\frac{2N^2+1}{6},$$

$$\Delta_{jj} = \frac{-\xi_j}{2(1-\xi_j^2)}, \; j=1,...,N-1,$$

$$\Delta_{ij} = \frac{\lambda_i}{\lambda_j}\frac{(-1)^{i+j}}{(\xi_i-\xi_j)}, i\neq j, i,j=1,...,N-1,$$

$$\lambda_i = \begin{cases} 2 & \text{if } i = 0, N \\ 1 & \text{otherwise} \end{cases}.$$

derived in [12].

We made use of the conformal transformation

$$r(\xi) = \frac{[1 + b\exp(-a)]r_{max}}{\left[1 + b\exp\left(-a\frac{1-\xi}{2}\right)\right]}\left(\frac{1-\xi}{2}\right) \quad (10)$$

that maps the standard interval $\xi \in [-1,1]$ onto the physical range of our problem $r \in [0, r_{max}]$. Because large matrices are involved, we numerically solved the eigenvalue problem using the Arnoldi type algorithm [11], which provides entire eigenvalue and eigenvector spectrum.

## 4 Modal Collocation With Orthogonal Basis For Inviscid Stability Analysis

The collocation method that we present in this section has the peculiar feature that can approximate the perturbation field for all types of boundary conditions, especially when the boundary limits are described by sophisticated expressions. We consider the mathematical model of an inviscid columnar vortex derived in [13] whose velocity profile is written as $\underline{V}(r) = [U(r), 0, W(r)]$.

$$G' + \frac{G}{r} + \frac{mH}{r} + kF = 0, \quad (11)$$

$$\left(\omega - \frac{mW}{r} - kU\right)G - \frac{2WH}{r} + P' = 0, \quad (12)$$

$$\left(-\omega + \frac{mW}{r} + kU\right)H + \left(W' + \frac{W}{r}\right)G + \frac{mP}{r} = 0, \quad (13)$$

$$\left(-\omega + \frac{mW}{r} + kU\right)F + U'G + kP = 0. \quad (14)$$

We assume for for this model that the radial amplitude of the velocity perturbation at the wall is negligible, i.e. $G(r_{max}) = 0$, for a truncated radius distance $r_{max}$ selected large enough such that the numerical results do not depend on that truncation of infinity. We have at $r = 0$

$$(|m| > 1), \quad F = G = H = P = 0, \quad (15)$$

$$(m = 0), \quad G = H = 0, F, P \text{ finite}, \quad (16)$$

$$(m = \pm 1), \quad H \pm G = 0, F = P = 0. \quad (17)$$

and at $r = r_{max}$

$$(|m| > 1), \quad F = G = H = P = 0, \quad (18)$$

$$(m = 0), \quad \frac{2W_{r_{max}}H}{r_{max}} - P' = 0, G = 0,$$

$$HkU_{r_{max}} - \omega H = 0, FkU_{r_{max}} - \omega F + kP = 0, \quad (19)$$

$$(m = \pm 1), \quad \frac{2W_{r_{max}}H}{r_{max}} - P' = 0, G = 0,$$

$$r_{max}H(kU_{r_{max}} - \omega) \pm HW_{r_{max}} \pm P = 0 = 0,$$

$$r_{max}F(kU_{r_{max}} - \omega) \pm FW_{r_{max}} + kr_{max}P = 0, \quad (20)$$

where $U_{r_{max}}$ and $W_{r_{max}}$ are the axial, respectively the tangential velocity calculated at domain limit $r_{max}$.

A different approach is obtained by taking as basis functions simple linear combinations of orthogonal polynomials. These are called bases of *modal* type, i. e., such that each basis function provides one particular pattern of oscillation of lower and higher frequency. We approximate the perturbation amplitudes as a truncated series of shifted Chebyshev polynomials

$$(F, G, H, P) = \sum_{k=1}^{N}(f_k, g_k, h_k, p_k) \cdot T_k^*, \quad (21)$$

where $T_k^*$ are shifted Chebyshev polynomials on the physical domain $[0, r_{max}]$. The clustered Chebyshev Gauss grid $\Xi = (\xi_j)_{1 \le j \le N}$ in $[-1, 1]$ is defined by relation

$$\xi_{j+1} = \cos\frac{\pi(j + N - 1)}{N - 1}, \quad \xi_{j+1} \in [-1, 1], \quad j = 0..N-1. \quad (22)$$

This formula has the advantage that in floating-point arithmetic it yields nodes that are perfectly symmetric about the origin, being clustered near the boundaries diminishing the negative effects of the Runge phenomena [11, 12]. This collocation nodes are the roots of Chebyshev polynomials and distribute the error evenly and exhibit rapid convergence rates with increasing numbers of terms.

In order to approximate the derivatives of the unknown functions, we express the derivative of the shifted Chebyshev polynomial $T_n^*$ as a difference between the previous and the following term

$$T_n^{*'}(r) = \frac{r_{max}}{4}\frac{(n-1)}{r(r_{max} - r)}\left[T_{n-1}^*(r) - T_{n+1}^*(r)\right], \quad n \ge 2. \quad (23)$$

Let us consider

$$F(r) = f_1 T_1^*(r) + \sum_{k=2}^{N} f_k T_k^*(r). \quad (24)$$

By differentiating (24) results

$$F'(r) = f_1 T_1^{*'}(r) + \sum_{k=2}^{N} f_k T_k^{*'}(r). \quad (25)$$

But $T_1^{*'}(r) = 0$ and involving relation (23) results

$$F'(r) = \sum_{k=2}^{N} f_k \frac{r_{max}}{4}\frac{(k-1)}{r(r_{max} - r)}\left[T_{k-1}^*(r) - T_{k+1}^*(r)\right]. \quad (26)$$

The interpolative differentiation matrix D that approximates the discrete derivatives has the elements

$$D_{m,n} = E_n(r_m), m = 2..N - 1, n = 2..N - 1, \quad (27)$$

where for $k = 2..N - 1$

$$E_k(r) = \frac{(k-1)}{r(r_{\max}-r)}\left[T^*_{k-1}(r) - T^*_{k+1}(r)\right]. \quad (28)$$

The eigenvalue problem governing the inviscid stability analysis appears now as a system of $4N$ equations, when include the boundary conditions. A special situation occur for the cases $m = \pm 1$, when only 7 relations define the boundary conditions. To regain the 8$^{th}$ equation we choose the third relation from the mathematical model and we compute it in the extreme node $r = r_{\max}$.

We have chosen this relation for few reasons. We observed that the equations that not contain the axial perturbation $F$ are the second and the third. The second equation contains the derivative of the pressure perturbation that cannot be computed in extreme nodes because the interpolative derivative matrix may produce singularities since contains the expression $E_k(r)$. The remain possibility is actually the third equation symmetrized.

The hydrodynamic model reads, for $j = 2..N-1$

$$G' + \frac{1}{r_j}\sum_{k=1}^{N} g_k T^*_k(r_j) + \frac{m}{r_j}\sum_{k=1}^{N} h_k T^*_k(r_j) + k\sum_{k=1}^{N} f_k T^*_k(r_j), \quad (29)$$

$$\left[\omega - \frac{mW}{r_j} - kU\right]\sum_{k=1}^{N} g_k T^*_k(r_j) - \frac{2W}{r_j}\sum_{k=1}^{N} h_k T^*_k(r_j) + P' = 0, (30)$$

$$\left[-\omega + \frac{mW}{r_j} + kU\right]\sum_{k=1}^{N} h_k T^*_k(r_j) +$$
$$+\left[W' + \frac{W}{r_j}\right]\sum_{k=1}^{N} g_k T^*_k(r_j) + \frac{m}{r_j}\sum_{k=1}^{N} p_k T^*_k(r_j) = 0, (31)$$

$$\left[-\omega + \frac{mW}{r_j} + kU\right]\sum_{k=1}^{N} f_k T^*_k(r_j) +$$
$$+U'\sum_{k=1}^{N} g_k T^*_k(r_j) + k\sum_{k=1}^{N} p_k T^*_k(r_j) = 0, \quad (32)$$

$$kr_{\max}U_{r\max}\sum_{k=1}^{N} h_k + (mW_{r\max} - r_{\max}\omega)\sum_{k=1}^{N} h_k +$$
$$+(W_{r\max} + r_{\max}W'_{r\max})\sum_{k=1}^{N} g_k + m\sum_{k=1}^{N} p_k = 0, \quad (33)$$

$$\sum_{1}^{N}(-1)^{k+1}g_k \pm \sum_{1}^{N}(-1)^{k+1}h_k = 0, \quad (34)$$

$$\sum_{1}^{N}(-1)^{k+1}f_k = \sum_{1}^{N}(-1)^{k+1}p_k = 0, \quad (35)$$

$$\frac{2W_{r\max}}{r_{\max}}\sum_{1}^{N} h_k - p_2 \frac{2}{r_{\max}} - \sum_{\substack{3 \\ k\ odd}}^{N} p_k \frac{2(k-1)}{r_{\max}}\left[\sum_{\substack{r=k-1 \\ k\ even}}^{2} 2\right] -$$
$$-\sum_{\substack{4 \\ k\ even}}^{N} p_k \frac{2(k-1)}{r_{\max}}\left[\sum_{\substack{r=k-1 \\ k\ odd}}^{2} 2 + 1\right] = 0, \quad (36)$$

$$\sum_{1}^{N} g_k = 0, \quad (37)$$

$$kU_{r\max}r_{\max}\sum_{1}^{N} h_k + (\pm W_{r\max} - \omega r_{\max})\sum_{1}^{N} h_k \pm \sum_{1}^{N} p_k = 0, \quad (38)$$

$$k\left(U_{r\max}r_{\max}\sum_{1}^{N} f_k + r_{\max}\sum_{1}^{N} p_k\right) + (\pm W_{r\max} - \omega r_{\max})\sum_{1}^{N} f_k = 0. (39)$$

Let us denote by $[r] = diag(r_i)$, $\left[\frac{1}{r}\right] = diag(1/r_i)$,
$[\eta] = (\eta_{ij})_{\substack{2\le i\le N-1 \\ 1\le j\le N}}$, $\eta_{ij} = T^*_j(r_i)$, $[U] = diag(U(r_i))$
$[W] = diag(W(r_i))$, $2 \le i \le N-1$. Written in matrix formulation, the hydrodynamic model reads

$$(kM_k + \omega M_\omega + mM_m + M_0)\bar{s} = 0,$$
$$\bar{s} = (f_1,...,f_N, g_1,...,g_N, h_1,...,h_N, p_1,...,p_N)^T, \quad (40)$$

where $M_k$, $M_\omega$, $M_m$ and $M_0$ are square matrices of dimension $4N$ and the elements being matrix blocks

$$M_k = \begin{pmatrix} \widetilde{M_k} \\ boundary\ conditions\ blocks \end{pmatrix},$$

$$M_\omega = \begin{pmatrix} \widetilde{M_\omega} \\ boundary\ conditions\ blocks \end{pmatrix},$$

$$M_m = \begin{pmatrix} \widetilde{M_m} \\ boundary\ conditions\ blocks \end{pmatrix},$$

$$M_0 = \begin{pmatrix} \widetilde{M_0} \\ boundary\ conditions\ blocks \end{pmatrix},$$

$$\widetilde{M_k} = \begin{pmatrix} [r][\eta] & 0 & 0 & 0 \\ 0 & [U][\eta] & 0 & 0 \\ 0 & 0 & [rU][\eta] & 0 \\ [U][\eta] & 0 & 0 & [\eta] \end{pmatrix}$$

$$\widetilde{M_\omega} = \begin{pmatrix} 0 & 0 & 0 & 0 \\ 0 & -[\eta] & 0 & 0 \\ 0 & 0 & -[r][\eta] & 0 \\ -[\eta] & 0 & 0 & 0 \end{pmatrix},$$

$$\widetilde{M_m} = \begin{pmatrix} 0 & 0 & [\eta] & 0 \\ 0 & \left[\frac{W}{r}\right][\eta] & 0 & 0 \\ 0 & 0 & [W][\eta] & [\eta] \\ \left[\frac{W}{r}\right][\eta] & 0 & 0 & 0 \end{pmatrix}$$

$$\widetilde{M_0} = \begin{pmatrix} 0 & [\eta] + [r]D & 0 & 0 \\ 0 & 0 & 2\left[\frac{W}{r}\right][\eta] & -D \\ 0 & [W][\eta] + [rW'][\eta] & 0 & 0 \\ 0 & [U'][\eta] & 0 & 0 \end{pmatrix}$$

where $D$ represents the interpolative derivative matrix.

## 5 Model Validation On a Q-Vortex Profile

Assuming the velocity profile of Q-Vortex, written in form

$$U(r) = a + e^{-r^2}, \quad W(r) = \frac{q}{r}\left(1 - e^{-r^2}\right), \quad (41)$$

where $q$ represents the swirl number and $a$ provides a measure of free-stream axial velocity, we perform a spatial stability analysis using the collocation method described above. The spectra of the eigenvalue problem governing the spatial stability is depicted in Figure 1.

It is noticeable that the eigenvalue with the largest imaginary part defines the most unstable mode. In Table 1 we have compared the results obtained by this method with those of Olendraru et al. [14], in the non axisymmetrical case $|m|>1$.

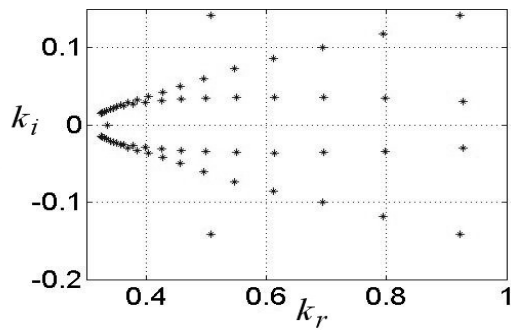

Fig.1 Spectra of the hydrodynamic eigenvalue problem computed at $\omega=0.01$, $m=-3$, $a=0$, $q=0.1$, for $N=100$ collocation nodes.

Table 1. Comparative results of the most amplified k-spatial wave at $a=0$, $q=0.1$, $\omega=0.01$ for the case of the counter-rotating mode $m=-3$: eigenvalue with largest imaginary part $k_{cr}=(k_r,k_i)$ and critical distance of the most amplified perturbation $r_c$.

| Shooting method [14] | |
|---|---|
| $k_{cr}=(0.506,-0.139)$ | $r_c=1.0005$ |
| Collocation method | |
| $k_{cr}=(0.50819,-0.14192)$ | $r_c=0.971$ |
| Error  0.79% | 2.94% |

## 6 Conclusion

In this paper we developed hydrodynamic models using spectral differential operators to investigate the spatial stability of swirling fluid systems, using two different methods.

When viscosity is considered as a valid parameter of the fluid, the hydrodynamic model is implemented using a nodal Lagrangean basis and the eigenvalue problem describing the viscous spatial stability is solved using the companion vector method. The second model for inviscid study is assessed for the construction of a certain class of shifted orthogonal expansion functions. The choice of the grid and the choice of the trial basis eliminate the singularities and the spectral differentiation matrix was derived to approximate the discrete derivatives. The models were applied to a Q-vortex structure, the scheme based on shifted Chebyshev polynomials providing good results.